





\input amstex
\documentstyle{amsppt}
\magnification=\magstep1
\parindent=1em
\baselineskip 15pt
\documentstyle{amsppt}
  \hsize=6.2truein
  \vsize=9.0truein
  \hoffset 0.1truein
  \parindent=2em
\NoRunningHeads
\pageno=1



\newcount\theTime
\newcount\theHour
\newcount\theMinute
\newcount\theMinuteTens
\newcount\theScratch
\theTime=\number\time
\theHour=\theTime
\divide\theHour by 60
\theScratch=\theHour
\multiply\theScratch by 60
\theMinute=\theTime
\advance\theMinute by -\theScratch
\theMinuteTens=\theMinute
\divide\theMinuteTens by 10
\theScratch=\theMinuteTens
\multiply\theScratch by 10
\advance\theMinute by -\theScratch

\def\today{{\number\day\space
 \ifcase\month\or
  January\or February\or March\or April\or May\or June\or
  July\or August\or September\or October\or November\or December\fi
 \space\number\year}}






















\define\biggnm#1{                            
  \bigg|\bigg|#1\bigg|\bigg|}

\define\bignm#1{                             
  \big|\big|#1\big|\big|}















\define\Cpx{\bold C}                         

\define\diag{\text{\rm diag}}                









\define\eqdef{{\;\overset\text{def}\to=\;}}     







\define\fpamalg#1{{\dsize\;                  
     \operatornamewithlimits*_{#1}\;}}


\define\freeprod#1#2{\mathchoice             
     {\operatornamewithlimits{\ast}_{#1}^{#2}}
     {\raise.5ex\hbox{$\dsize\operatornamewithlimits{\ast}
      _{#1}^{#2}$}\,}
     {\text{oops!}}{\text{oops!}}}

\define\freeprodi{\mathchoice                
     {\operatornamewithlimits{\ast}
      _{\iota\in I}}
     {\raise.5ex\hbox{$\dsize\operatornamewithlimits{\ast}
      _{\sssize\iota\in I}$}\,}
     {\text{oops!}}{\text{oops!}}}













\define\Hil{{\mathchoice                     
     {\text{\eusm H}}
     {\text{\eusm H}}
     {\text{\eusms H}}
     {\text{\eusmss H}}}}


















\define\Jc{{\Cal J}}                         








\define\lrnm#1{\left|\left|#1\right|\right|} 

\define\lspan{\text{\rm span}@,@,@,}         










\define\nm#1{||#1||}                         


\define\Nats{\Naturals}

\define\Naturals{{\bold N}}                  








\define\owedge{{                             
     \operatorname{\raise.5ex\hbox{\text{$
     \ssize{\,\bigcirc\llap{$\ssize\wedge\,$}\,}$}}}}}

\define\owedgeo#1{{                          
     \underset{\raise.5ex\hbox
     {\text{$\ssize#1$}}}\to\owedge}}










\define\QED{\newline                         
            \line{$\hfill$\qed}\enddemo}











\define\restrict{\lower .3ex                 
     \hbox{\text{$|$}}}




\def\sigmabar{{\overline\sigma}}

\define\smd#1#2{\underset{#2}\to{#1}}          

\define\smdb#1#2{\undersetbrace{#2}\to{#1}}    

\define\smdbp#1#2#3{\overset{#3}\to            
     {\smd{#1}{#2}}}

\define\smdbpb#1#2#3{\oversetbrace{#3}\to      
     {\smdb{#1}{#2}}}

\define\smdp#1#2#3{\overset{#3}\to             
     {\smd{#1}{#2}}}

\define\smdpb#1#2#3{\oversetbrace{#3}\to       
     {\smd{#1}{#2}}}

\define\smp#1#2{\overset{#2}\to                
     {#1}}                                     










\define\tocdots                              
  {\leaders\hbox to 1em{\hss.\hss}\hfill}    














\def\wbar{{\overline w}}                     

















  \newcount\mycitestyle \mycitestyle=1 

  \newcount\bibno \bibno=0
  \def\newbib#1{\advance\bibno by 1 \edef#1{\number\bibno}}
  \ifnum\mycitestyle=1 \def\cite#1{{\rm[\bf #1\rm]}} \fi
  \def\scite#1#2{{\rm[\bf #1\rm, #2]}}


  \newcount\ignorsec \ignorsec=0

  \newcount\secno \secno=0
  \def\newsec#1{\procno=0 \subsecno=0 \ignorsec=0
    \advance\secno by 1 \edef#1{\number\secno}
    \edef\currentsec{\number\secno}}

  \newcount\subsecno
  \def\newsubsec#1{\procno=0 \advance\subsecno by 1 \edef#1{\number\subsecno}
    \edef\currentsec{\number\secno.\number\subsecno}}

  \newcount\appendixno \appendixno=0
  \def\newappendix#1{\procno=0 \ignorsec=0 \advance\appendixno by 1
    \ifnum\appendixno=1 \edef\appendixalpha{\hbox{A}}
      \else \ifnum\appendixno=2 \edef\appendixalpha{\hbox{B}} \fi
      \else \ifnum\appendixno=3 \edef\appendixalpha{\hbox{C}} \fi
      \else \ifnum\appendixno=4 \edef\appendixalpha{\hbox{D}} \fi
      \else \ifnum\appendixno=5 \edef\appendixalpha{\hbox{E}} \fi
      \else \ifnum\appendixno=6 \edef\appendixalpha{\hbox{F}} \fi
    \fi
    \edef#1{\appendixalpha}
    \edef\currentsec{\appendixalpha}}

  \newcount\procno \procno=0
  \def\newproc#1{\advance\procno by 1
   \ifnum\ignorsec=0 \edef#1{\currentsec.\number\procno}
                     \edef\currentproc{\currentsec.\number\procno}
   \else \edef#1{\number\procno}
         \edef\currentproc{\number\procno}
   \fi}

  \newcount\subprocno \subprocno=0
  \def\newsubproc#1{\advance\subprocno by 1
   \ifnum\subprocno=1 \edef#1{\currentproc a} \fi
   \ifnum\subprocno=2 \edef#1{\currentproc b} \fi
   \ifnum\subprocno=3 \edef#1{\currentproc c} \fi
   \ifnum\subprocno=4 \edef#1{\currentproc d} \fi
   \ifnum\subprocno=5 \edef#1{\currentproc e} \fi
   \ifnum\subprocno=6 \edef#1{\currentproc f} \fi
   \ifnum\subprocno=7 \edef#1{\currentproc g} \fi
   \ifnum\subprocno=8 \edef#1{\currentproc h} \fi
   \ifnum\subprocno=9 \edef#1{\currentproc i} \fi
   \ifnum\subprocno>9 \edef#1{TOO MANY SUBPROCS} \fi
  }

  \newcount\tagno \tagno=0
  \def\newtag#1{\advance\tagno by 1 \edef#1{\number\tagno}}




\def\dim{\text{dim }}

\def\diag{\text{diag }}



\def\Com{\text{\rm Com }}
\redefine\Hil{{\Cal H}}
\redefine\Naturals{{\Bbb N}}


\newsec{\SingleComm}
 \newproc{\UpperTriang}
 \newproc{\InIdeal}
 \newproc{\SingleNec}
 \newproc{\UppTrCom}
 \newproc{\QnilEx}
  \newtag{\prodvjuj}
  \newtag{\thetasinth}
\newsec{\SpectralTraces}
 \newproc{\SpecOut}
  \newtag{\specout}
  \newtag{\etalambda}
 \newproc{\WhatsLeftNil}
 \newproc{\QuasiNTriang}
 \newproc{\DiagQnil}
 \newproc{\GSDiagQnil}
 \newproc{\SpecTrace}

\topmatter

\title Spectral characterization of sums of commutators II \endtitle

\author
  K.J. Dykema and N.J.Kalton
    \footnote{\hbox{Supported by NSF Grant DMS--9500125\hfil}}
\endauthor

\date 27 August, 1997 \enddate

\address1
Department of Mathematics and Computer Science,
Odense Universitet,
Campusvej 55
DK-5230 Odense M,
Denmark
\endaddress
\email  dykema\@imada.ou.dk
\endemail

\address2
Department of Mathematics,
University of Missouri,
Columbia, Mo.  65211, U.S.A.
\endaddress
\email  nigel\@math.missouri.edu
\endemail


\abstract
  For countably generated ideals, $\Jc$, of $B(\Hil)$, geometric stability is
  necessary for the canonical spectral characterization of sums
  of $(\Jc,B(\Hil))$--commutators to hold.
  This answers a question raised by
  Dykema, Figiel, Weiss and Wodzicki.
  There are some ideals, $\Jc$, having quasi--nilpotent elements that are not
  sums of $(\Jc,B(\Hil))$--commutators.
  Also, every trace on every geometrically stable ideal is a spectral
  trace.
\endabstract

\endtopmatter

\document \TagsOnRight \baselineskip=18pt

\heading Introduction \endheading

Let $\Cal H$ be a separable infinite--dimensional Hilbert space, and let
$\Cal J$ be a (two-sided) ideal contained in the ideal of compact
operators
$\Cal K(\Cal H)$ on $\Cal H.$
We define the {\it commutator subspace} $\Com \Cal J$ to be the
linear span of commutators $[A,B]=AB-BA$ where $A\in\Cal J$ and $B\in\Cal
B(\Cal H).$

In the immediately preceding paper,~\cite{5} N.J.~Kalton showed that, for a
wide class of ideals including quasi--Banach ideals (i.e. ideals equipped
with a complete ideal quasi--norm) $\Com\Jc$ has the following spectral
characterization.
For
$T\in\Cal J$, let $\lambda_k=\lambda_k(T)$ be the eigenvalues of $T$ listed
according to algebraic multiplicity and such that
$|\lambda_1|\ge|\lambda_2|\ge\cdots$.
Then $T\in\Com\Cal J$ if and only if
$\diag\{\frac1n\sum_{k=1}^n\lambda_k\}\in\Cal J.$
The sufficient condition from~\cite{5} for $\Com\Cal J$ to have a spectral
characterization in the above sense is {\it
geometric stability}, i.e. the condition that if $(s_n)$ is a monotone
decreasing real sequence then $\diag\{s_n\}\in\Cal J$ if and only if
$\diag\{(s_1\ldots s_n)^{1/n}\}\in\Cal J$.
(See the introduction of~\cite{5} for background.)

In~\S\SingleComm{} of this paper, we show that when $\Jc$ is a countably generated
ideal, then
geometric stability is a necessary condition for this spectral
characterization of $\Com\Jc$.
In particular, for every countably generated ideal, $\Jc$, which is not
geometrically stable, we find an element $T$, of $\Com\Jc$ for which
$\diag\{\frac1n(\lambda_1(T)+\cdots+\lambda_n(T))\}\not\in\Jc$.
We also give examples of ideals for which
spectral characterization of $\Com\Jc$ fails in the opposite direction, and in
in an extreme way, in that $\Jc$ has a quasi--nilpotent element that is not in
$\Com\Jc$.

In~\S\SpectralTraces, we show, based on an elementary decomposition
result for compact operators,
that every trace, $\tau$, on an arbitrary  geometrically stable ideal, $\Jc$,
is a spectral trace, i.e\. for every $T\in\Jc$,
$\tau(T)$ depends only on the eigenvalues of $T$ and their multiplicities.

\heading \SingleComm.  Countably generated ideals \endheading

\proclaim{Proposition \UpperTriang}Suppose $s_1\ge s_2\ge\cdots\ge 0.$  Suppose
that
$(\lambda_n)_{n=1}^{\infty}$ is any sequence of complex numbers with
$|\lambda_1|\ge |\lambda_2|\ge \cdots$ and such that
$|\lambda_1\ldots\lambda_n|\le s_1\ldots s_n$ for all $n.$  Then there is
an upper triangular operator $A=(a_{jk})_{j,k}$ such that $s_n(A)\le
s_n$ for all $n$ and
$a_{jj}=\lambda_j$ for all $j.$\endproclaim

\demo{Proof}This is the infinite dimensional extension of a result of
Horn \cite{3} (see Gohberg--Krein \cite{2} Remark II.3.1).  By Horn's result we
can find for each $n$ an upper--triangular matrix
$A^{(n)}=(a^{(n)}_{jk})$ such that
$a^{(n)}_{jk}=0$ if $\max{j,k}>n$, $a^{n}_{jj}=\lambda_j$ for $1\le j\le
n$ and $s_j(A^{(n)})\le s_j$ for $1\le j<\infty.$

Now we can pass to a subsequence $(B^{(n)})$ of $(A^{(n)})$ which is
convergent to the weak operator topology to some operator $A=(a_{ij})$
which is upper--triangular and has $a_{jj}=\lambda_j$ for all $j.$  It
remains to show that $s_j(A)\le s_j$ for all  $j.$

        Fix $j>1$ (the case $j=1$ is well known.)  For each $n$ there is a
subspace
$E_n$ of dimension
$j-1$
so that if $P_n$ is the orthogonal projection with kernel $E_n$ then
$\|B^{(n)}P_n\|\le s_j.$  Let $\{f^n_k\}_{k=1}^{j-1}$ be an orthonormal
basis of $E_n.$  By passing to a further subsequence we can suppose that
$\lim_{n\to\infty}f^n_k=f_k$ exists weakly for each $1\le k\le j-1.$
Let $E$ be the span of $(f_k)_{k=1}^{j-1}$, so that $\dim E\le
j-1.$ Suppose $x\in E^{\perp}.$  Then $\lim (x,f^n_k)=0$ for $1\le k\le
j-1.$  Hence $\lim_{n\to\infty}\|x-P_nx\|=0$.  Thus $\|Ax\|\le
\liminf_{n\to\infty}\|B^{(n)}x\|  \le \liminf_{n\to\infty}\|B^{(n)}P_nx\|\le
s_j\|x\|.$   This shows that $s_j(A)\le s_j.$\qed\enddemo

Let us introduce a bit of notation.
Suppose $u=(u_j)_{j=1}^\infty$ and $t=(t_j)_{j=1}^\infty$ are
decreasing sequences of positive numbers.
We write $t\le u$ to mean that $u_j\le t_j$ for every $j$.
For $c>0$ we let $cu$ denote the sequence $(cu_j)_{j=1}^\infty$.
We let
$u\oplus u$ denote the sequence $(u_1,u_1,u_2,u_2,u_3,u_3,\ldots)$.

\proclaim{Lemma \InIdeal}Let $u=(u_j)_{j=1}^\infty$ and $t=(t_j)_{j=1}^\infty$ be
sequences of positive numbers decreasing to zero such that $\diag\{t_j\}$ is in
the ideal of compact operators generated by $\diag\{u_j\}$.
Then $\diag\{(t_1t_2\cdots t_k)^{1/k}\}$ is in the ideal of compact operators
generated by $\diag\{(u_1u_2\ldots u_k)^{1/k}\}$.
\endproclaim

\demo{Proof}The hypotheses imply that there is a positive integer $n$ and a
positive number $c$ such that $t_{nj}\le cu_j$ for every $j$, and we may suppose
$n$ is a power of $2$.
Hence, in order to prove the lemma it will suffice to prove it in each of the
following special cases, in turn.
\roster
\item"(i)" $t\le u$,
\item"(ii)" $t=cu$ for $c>1$,
\item"(iii)" $t=u\oplus u$.
\endroster

The cases (i) and (ii) are clear.
If (iii) holds then
$(t_1\ldots t_{2n})^{1/2n}=(u_1\ldots u_n)^{1/n}$,
which shows that
$\diag\{(t_1\ldots t_k)^{1/k}\}$ is in the ideal generated by
$\diag\{(u_1\ldots u_n)^{1/n}\}$, and the lemma is proved.
\QED

 \proclaim{Theorem \SingleNec}Let $\Cal J$ be a countably generated ideal in
$\Cal K(\Cal H).$  Then the following conditions are equivalent:
\roster
\item"(i)" $\Cal J$ is geometrically stable.
\item"(ii)" If $T\in\Cal J$ then $\diag\{\lambda_n\}\in\Cal J$ where
$\lambda_n=\lambda_n(T)$.
\item"(iii)" If $T\in\Com\Cal J$ then $\diag
\{\frac1n(\lambda_1+\cdots+\lambda_n)\}\in\Cal J,$ for some ordering of
$\lambda_n$ so that $(|\lambda_n|)$ is decreasing.
\endroster
\endproclaim

\demo{Proof} (i) implies (iii) is proved for any ideal in \cite{5}.
That (i) implies (ii) is also noted in \cite{5}; this follows from
the inequalities $|\lambda_n(T)|\le (s_1(T)\ldots s_n(T))^{1/n}.$

Since $\Cal J$ is countably generated there is a countable family,
$\{u^{(k)}\mid k\in\Naturals\}$, of decreasing sequences of positive numbers
such that $T\in\Jc$ if and only if for some $k$ we have $s_n(T)\le u^{(k)}_n$
for all $n\ge1$.
Indeed, let
$\{u^{(k)}\mid k\in\Naturals\}=\{w^{(i,j)}\mid i,j\in\Naturals\}$
where $w^{(i,j)}_n=jw^{(i)}_{n/j}$ and where
$w^{(i)}_n=s_n(T_1)+s_n(T_2)+\cdots+s_n(T_i)$, with
$\{T_i\mid i\in\Naturals\}$ a countable generating set of the ideal $\Jc$.

We will prove that each of  (ii) and
(iii) individually implies
(i) by supposing that (i) fails and showing that both (ii) and (iii) fail.
If (i) fails then there is a sequence $(t_n)_{n=1}^{\infty}$
such that $\diag\{t_n\}\in\Jc$ but
$\diag\{(t_1t_2\cdots t_n)^{1/n}\}\notin\Jc$.
Now for every $k\in\Naturals$, since $\diag\{(t_{n+k})_{n=1}^\infty\}$ generates
the same ideal as $\diag\{t_n\}$, it follows from Lemma~\InIdeal{} that
$\diag\bigl\{\bigl((t_{k+1}t_{k+2}\cdots
t_{k+n})^{1/n}\bigr)_{n=1}^\infty\bigr\}\notin\Jc$.

It follows that we can find a sequence $m_n\uparrow \infty$ such that
if $m_0=0$ and $p_n=m_n-m_{n-1}$ for $n\ge 1$ then $p_n/m_n\to 1$ and
$$ \wbar_n\eqdef\left(\prod_{k=m_{n-1}+1}^{m_n}t_k\right)^{1/p_n}>
u^{(n)}_{m_n}.$$
Let $w_k=\wbar_n$ for every
$m_{n-1}+1\le k\le m_n.$
We may assume without loss of generality that $\wbar_n>\wbar_{n+1}$ for every
$n\in\Bbb N$.
By Proposition~\UpperTriang{} there is an upper triangular
operator $A\in\Cal J$ with diagonal entries $a_{kk}=w_k$ for $k\in \Bbb
N.$   Now $\diag (w_k)\notin \Cal J$ by construction and
hence (ii) fails.
Therefore (ii) implies (i).

We now show that also (iii) fails, which will finish the proof of the theorem.
We can choose $\sigmabar_1>\sigmabar_2>\cdots>0$ decreasing quickly enough so
that $\wbar_n-\sigmabar_n>\wbar_{n+1}$ for every $n$ and such that, letting
$\sigma_k=\sigmabar_n$ for every $m_{n-1}+1\le k\le m_n$, we have
$\diag(\sigma_k)\in\Cal J$.
Define for each $n\in\Bbb N$
$$ \epsilon_{2n-1}=\epsilon_{2n}={\tsize\frac12}
\min(\sigmabar_{2n-1}-\sigmabar_{2n},\sigmabar_{2n}-\sigmabar_{2n+1}). $$
Then for every $n\in\Bbb N$ we have that $\epsilon_n>0$, that
$\sigmabar_n-\epsilon_n>\sigmabar_{n+1}$ and that
$$ \left|\sum_{j=1}^n(-1)^{j-1}\epsilon_j\right|\le
\cases\sigmabar_n&\text{if $n$ is odd}\\
0&\text{if $n$ is even}.\endcases $$
Let $(\mu_k)_{k=1}^\infty$ be the sequence defined by
$$ \mu_k=\cases
\sigmabar_n&\text{if }m_{n-1}+1\le k\le m_n\text{ and $n$ is odd}\\
\sigmabar_n-(\epsilon_n/p_n)&\text{if }m_{n-1}+1\le k\le m_n
\text{ and $n$ is even}\\ \endcases $$
and let  $(\nu_k)_{k=1}^\infty$ be the sequence defined by
$$ \nu_k=\cases
\sigmabar_n-(\epsilon_n/p_n)&\text{if }m_{n-1}+1\le k\le m_n
\text{ and $n$ is odd}\\
\sigmabar_n&\text{if }m_{n-1}+1\le k\le m_n\text{ and $n$ is even.}
\endcases $$
Note that $(u_k)_{k=1}^\infty$ and $(\nu_k)_{k=1}^\infty$
are decreasing sequences.
Let $D_1=\diag\{\mu_k\}$ and $D_2=\diag\{\nu_k\}$.
We will show that $D_1\oplus(-D_2)\in\Com\Cal J$.
Let $(\lambda'_k)_{k=1}^\infty$ be the eigenvalue sequence, arranged in some
order of decreasing absolute value, of the operator $D_1\oplus(-D_2)$.
We will show that $\diag\{\frac1k(\lambda'_1+\cdots+\lambda'_k)\}\in\Cal J$,
which by~\cite{1} shows $D_1\oplus(-D_2)\in\Com\Cal J$.
The sequence $(\lambda'_k)_{k=1}^\infty$ is
$$ \align
\undersetbrace{p_1\text{ times}}\to{\sigmabar_1,\ldots,\sigmabar_1},
&\undersetbrace{p_1\text{ times}}
 \to{-\sigmabar_1+{\tsize\frac{\epsilon_1}{p_1}},
     \ldots,-\sigmabar_1+{\tsize\frac{\epsilon_1}{p_1}}},\;
\undersetbrace{p_2\text{ times}}\to{-\sigmabar_2,\ldots,-\sigmabar_2},\;
\undersetbrace{p_2\text{ times}}
 \to{\sigmabar_2-{\tsize\frac{\epsilon_2}{p_2}},
     \ldots,\sigmabar_2-{\tsize\frac{\epsilon_2}{p_2}}},\;\ldots\\\vspace{2ex}
\ldots,\;&\undersetbrace{p_{2n-1}\text{ times}}
 \to{\sigmabar_{2n-1},\ldots,\sigmabar_{2n-1}},\;
\undersetbrace{p_{2n-1}\text{ times}}
 \to{-\sigmabar_{2n-1}+{\tsize\frac{\epsilon_{2n-1}}{p_{2n-1}}},
     \ldots,-\sigmabar_{2n-1}+{\tsize\frac{\epsilon_{2n-1}}{p_{2n-1}}}}, \\\vspace{2ex}
&\undersetbrace{p_{2n}\text{ times}}\to{-\sigmabar_{2n},\ldots,-\sigmabar_{2n}},\;
\undersetbrace{p_{2n}\text{ times}}
 \to{\sigmabar_{2n}-{\tsize\frac{\epsilon_{2n}}{p_{2n}}},
     \ldots,\sigmabar_{2n}-{\tsize\frac{\epsilon_{2n}}{p_{2n}}}},\;\ldots.
\endalign $$
Let $\eta'_{k}=\sum_{j=1}^k\lambda'_k$.
Then for every $n\in\Bbb N$,
$$ \alignat2
\eta'_{2m_{2n-2}+k}&=k\sigmabar_{2n-1}&&(1\le k\le p_{2n-1}) \\\vspace{2ex}
\eta'_{2m_{2n-2}+p_{2n-1}+k}&=(p_{2n-1}-k)\sigmabar_{2n-1}
 +\frac{k\epsilon_{2n-1}}{p_{2n-1}}&&(1\le k\le p_{2n-1}) \\\vspace{2ex}
\eta'_{2m_{2n-1}+k}&=\epsilon_{2n-1}-k\sigmabar_{2n}&&(1\le k\le p_{2n}) \\
  \vspace{2ex}
\eta'_{2m_{2n-1}+p_{2n}+k}&=\epsilon_{2n-1}-(p_{2n}-k)\sigmabar_{2n}
 -\frac{k\epsilon_{2n}}{p_{2n}}\quad&&(1\le k\le p_{2n}).
\endalignat $$
Making crude (but sufficient) estimates, we get
$$ \alignat2
\left|\frac{\eta'_{2m_{2n-2}+k}}{2m_{2n-2}+k}\right|
 &\le\sigmabar_{2n-1}&&(1\le k\le p_{2n-1}) \\\vspace{2ex}
\left|\frac{\eta'_{2m_{2n-2}+p_{2n-1}+k}}{2m_{2n-2}+p_{2n-1}+k}\right|
 &\le2\sigmabar_{2n-1}&&(1\le k\le p_{2n-1}) \\\vspace{2ex}
\left|\frac{\eta'_{2m_{2n-1}+k}}{2m_{2n-1}+k}\right|
 &\le\sigmabar_{2n-1}+\sigmabar_{2n}&&(1\le k\le p_{2n}) \\\vspace{2ex}
\left|\frac{\eta'_{2m_{2n-1}+p_{2n}+k}}{2m_{2n-1}+p_{2n}+k}\right|
 &\le\sigmabar_{2n-1}+\sigmabar_{2n}.\quad&&(1\le k\le p_{2n})
\endalignat $$
Using $\diag\{\sigma_k\}\in\Cal J$ we see that
$\diag\{\frac1k(\lambda'_1+\cdots+\lambda'_k)\}\in\Cal J$.
Hence by~\cite{1} $D_1\oplus(-D_2)\in\Com\Cal J$.

Let $A\in\Cal J$ be the upper triangular operator with diagonal entries
$a_{kk}=w_k$ as constructed above.
Consider the operator
$$ T=(A+D_1)\oplus(-A-D_2). $$
Then $T\in\Com\Cal J$ because $A\oplus(-A)\in\Com\Cal J$ and
$D_1\oplus(-D_2)\in\Com\Cal J$.
Let $\lambda_k=\lambda_k(T)$ be the eigenvalues of $T$ listed according to
algebraic multiplicity and in order of decreasing absolute value.
Then
$$ \alignat2
\lambda_{2m_{2n-2}+k}&=\lambda'_{2m_{2n-2}+k}+\wbar_{2n-1}
 &&(1\le k\le p_{2n-1}) \\\vspace{1ex}
\lambda_{2m_{2n-2}+p_{2n-1}+k}&=\lambda'_{2m_{2n-2}+p_{2n-1}+k}-\wbar_{2n-1}
 \quad&&(1\le k\le p_{2n-1}) \\\vspace{1ex}
\lambda_{2m_{2n-1}+k}&=\lambda'_{2m_{2n-1}+k}-\wbar_{2n}
 &&(1\le k\le p_{2n}) \\\vspace{1ex}
\lambda_{2m_{2n-1}+p_{2n}+k}&=\lambda'_{2m_{2n-1}+p_{2n}+k}+\wbar_{2n}
 \quad&&(1\le k\le p_{2n}).
\endalignat $$
Let $\eta_{k}=\sum_{j=1}^k\lambda_j$.
Then
$$ \alignat2
\eta_{2m_{2n-2}+k}&=\eta'_{2m_{2n-2}+k}+k\wbar_{2n-1}
 &&(1\le k\le p_{2n-1}) \\\vspace{1ex}
\eta_{2m_{2n-2}+p_{2n-1}+k}
  &=\eta'_{2m_{2n-2}+p_{2n-1}+k}+(p_{2n-1}-k)\wbar_{2n-1}
 \quad&&(1\le k\le p_{2n-1}) \\\vspace{1ex}
\eta_{2m_{2n-1}+k}&=\eta'_{2m_{2n-1}+k}-k\wbar_{2n}
 &&(1\le k\le p_{2n}) \\\vspace{1ex}
\eta_{2m_{2n-1}+p_{2n}+k}&=\eta'_{2m_{2n-1}+p_{2n}+k}-(p_{2n}-k)\wbar_{2n}
 \quad&&(1\le k\le p_{2n}).
\endalignat $$

We will show that $\diag\{\frac1k(\lambda_1+\cdots+\lambda_k)\}\notin\Cal J$.
Since $\diag\{\frac1k\eta'_{k}\}\in\Cal J$, it will suffice to show that
$\diag\{\frac1k(\eta_{k}-\eta'_{k})\}\notin\Cal J$.
However, taking the absolute value of a subsequence of
$(\frac1k(\eta_{k}-\eta'_{k}))_{k=1}^\infty$ and
using $\lim_{n\to\infty}(p_n/m_n)=1$, we see that
$\diag\{\frac1k(\eta_{k}-\eta'_{k})\}\in\Cal J$ would imply
$\diag\{w_k\}\in\Cal J$, which would be a contradiction.
Therefore
$$ \diag\{\frac1k(\lambda_1+\cdots+\lambda_k)\}\notin\Cal J $$
and hence (iii) fails.
\QED

We will conclude by showing that it is in general impossible to
characterize membership in $\Com J$ by considering only the eigenvalues.
To this end we will construct a quasi--nilpotent operator $T$ so
that $T\notin \Com\Cal J_T.$

We will start by considering any singly generated ideal, $\Jc$.  Then there is
decreasing sequence of positive numbers $(u_n)$ such that
$T\in\Cal J$ if and only for some $C>0$ and $0<\alpha\le 1$ we have
$s_n(T)\le Cu_{\alpha n}.$

\proclaim{Lemma \UppTrCom}
Let $\Jc$ and $(u_n)$ be as above and
suppose $(v_n)_{n=1}^{\infty}$ and
$(w_n)_{n=1}^{\infty}$ are a pair of real sequences such that
\roster
\item"(i)" we have $v_1\ge v_2\ge \cdots\ge 0$;
\item"(ii)" for each $k\ge 1$ the sequence $(v_n)$ is constant on the dyadic
block
$$ \{2^{k-1},2^{k-1}+1,\ldots,2^k-1\};$$
\item"(iii)" for each $n\in \Bbb N$ we have $\prod_{k=1}^nv_k\le
\prod_{k=1}^nu_k$;
\item"(iv)" for each $n$ we have $|w_n|\le |u_n|$;
\item"(v)" for each $n$ we have $|\sum_{k=1}^nw_k|\le nv_n$.
\endroster
Then there is an upper-triangular operator $A\in\Com\Cal J$  with
$a_{nn}=w_n$ for all $n.$\endproclaim
\demo{Proof}The proof is essentially the same as for the corresponding
result for diagonal operators in \cite{4} and \cite{1}. We start by
setting
$$\eta_n=2^{1-k}\sum_{j=2^{k-1}}^{2^k-1}w_j$$
when $1\le k< \infty$ and $2^{k-1}\le n\le 2^k-1.$
Notice that
$$ |\eta_n-w_n| \le u_{n/2}$$ and that, if $2^{k-1}\le n\le 2^k-1,$
$$ |\sum_{j=1}^n(\eta_j-w_j)| =|\sum_{j=2^{k-1}}^n (\eta_j-w_j)|\le
2^{k}u_{n/2}$$
from which it follows that $\diag\{(\eta_n-w_n)\}\in\Com \Cal J.$

It therefore will suffice to show that there is an upper-triangular
operator $A$ in $\Com\Cal J$ with $a_{nn}=\eta_n$ for all $n.$
To this end we define a sequence $\xi_n$ by setting
$$ \xi_n =
2^{1-k}\sum_{j=1}^{2^k-1}\eta_j=2^{1-k}\sum_{j=1}^{2^k-1}w_j$$
for
$2^{k-1}\le n<2^{k-1}.$ By hypothesis we have $|\xi_n|\le 2v_n$ for every
$n\in \Bbb N.$  Now we note that by Proposition~\UpperTriang, there is an
upper-triangular matrix $C\in\Cal J$ with $c_{nn}=v_n$ for every $n$;
hence there is an  upper-triangular matrix $B\in\Cal J$ with
$b_{nn}=\xi_n$ for $n\in \Bbb N.$

Now consider the  isometries $U_1,U_2$ defined by $U_1e_n=e_{2n+1}$ and
$U_2e_n=e_{2n}.$
Consider the commutator $[U_1B,U_1^*].$  We have
$$ [U_1^*,U_1B](e_{2n+1}) = Be_{2n+1}-U_1Be_n,$$
$$ [U_1^*,U_1B](e_{2n})= Be_{2n}$$
and
$$ [U_1^*,U_1B]e_1= Be_1.$$
Similarly
$$ [U_2^*,U_2B]e_{2n}= Be_{2n}-U_2Be_n$$
$$ [U_2^*,U_2B]e_{2n+1}=Be_{2n+1}$$
and
$$ [U_2^*,U_2B]e_1=Be_1.$$

Thus if $A=\frac12([U_1^*,U_1B]+[U_2^*,U_2B])$ then $Ae_1=Be_1$ and
$$ \align
(Ae_{2n+1},e_{2m+1})&= (Be_{2n+1},e_{2m+1})-\frac12(Be_n,e_m) \\
(Ae_{2n+1},e_{2m}) &= (Be_{2n+1},e_{2m+1}\\
(Ae_{2n+1},e_1) &= (Be_{2n+1},e_1)=0 \\
(Ae_{2n},e_{2m+1})&= (Be_{2n},e_{2m+1})\\
(Ae_{2n},e_{2m})&= (Be_{2n},e_{2m}) -\frac12(Be_n,e_m)\\
(Ae_{2n},e_1)&=(Be_{2n},e_1)=0.
\endalign
$$

Thus $(Ae_n,e_m)=0$ if $m<n$ and $a_{11}=\xi_1=\eta_1$ while
if $2^{k-1}\le n\le 2^k-1$ with $k\ge 2$ we have
$a_{nn}=\xi_{2^{k-1}}-\frac12\xi_{2^{k-2}} = \eta_n.$  Since
$A\in\Com\Cal J,$ this completes the proof.\qed\enddemo

We now turn to the construction of some examples of ideals, $\Jc$, having
quasi-nilpotents that are not sums of commutators.

\proclaim{Examples \QnilEx}
Let $0=p_0<p_1<p_2<\cdots$ be integers such that
$p_{n+1}>p_n+2n2^{2n}$ for every $n\ge0$.
Let
$$ u_k=2^{-p_n}\qquad
\text{if}\quad2^{p_{n-1}}\le k\le 2^{p_n}-1,\quad(n\ge1) $$
and let $\Jc$ be the ideal of $B(\Hil)$ generated by
$\diag\{u_k\}_{k=1}^\infty$.
Then there is a quasi--nilpotent operator, $T\in\Jc$ such that
$T\notin\Com\Jc$.
\endproclaim
\demo{Proof}
Let $q_n=2n+p_n$ and
$$ v_k=2^n2^{-p_{n+1}}
\qquad\text{if}\quad2^{q_{n-1}}\le k\le 2^{q_n}-1,\quad(n\ge1). $$
We claim that
$$ \prod_{j=1}^kv_j\le\prod_{j=1}^ku_j \tag{\prodvjuj} $$
for every $k\in\Nats$.
Let $\log$ denote the base~$2$ logarithm.
For every $n\in\Nats$, we have
$$ \align
\log\left(\prod_{j=2^{q_{(n-1)}}}^{(2^{q_n})-1}v_ju_j^{-1}\right)=&
2^{p_n}\bigl((1-2^{-p_n+p_{n+1}+2n-2})(-p_{n+1}+p_n)+ \\
&+(1-2^{-p_n+p_{n-1}-2})n2^{2n}\bigr).
\endalign $$
But
$$ 1-2^{-p_n+p_{n-1}-2}<2(1-2^{-p_n+p_{n-1}+2n-2}) $$
so
$$ \log\left(\prod_{j=2^{q_{(n-1)}}}^{(2^{q_n})-1}v_ju_j^{-1}\right)
<2^{p_n}(1-2^{-p_n+p_{n-1}+2n-2})(-p_{n+1}+p_n+2n2^{2n})<0. $$
Therefore, by induction on $n$,~(\prodvjuj) holds whenever $k=2^{q_n}-1$.
But since $v_j<u_j$ when $2^{q_n}\le j\le 2^{p_{n+1}}-1$ and
$v_j>u_j$ when $2^{p_{n+1}}\le j\le 2^{q_{n+1}}-1$, it follows that~(\prodvjuj)
holds for all $k\in\Nats$.

Let $\sigma_k=\sum_{j=1}^ku_j$ and
$$ \theta_k=\inf_{j\in\Naturals}(jv_j+|\sigma_k-\sigma_j|). $$
Then $\theta_1\le v_1\le u_1$ and for $k\ge2$,
$|\theta_k-\theta_{k-1}|\le u_k$.
Let $w_1=\theta_1$ and $w_k=\theta_k-\theta_{k-1}$ ($k\ge2$).
Then for every $k\in\Nats$ we have $|w_k|\le u_k$ and
$0\le\sum_{j=1}^kw_j=\theta_k\le kv_k$.
Therefore, by Lemma~\UppTrCom{} there is an upper triangular $A\in\Com\Jc$ with
diagonal elements $a_{kk}=w_k$.
However, also $W\eqdef\diag\{w_k\}\in\Jc$.
Let $T=A-W$.
Then $T\in\Jc$ is quasi--nilpotent.
We will show that $T\notin\Com\Jc$ by showing $W\notin\Com\Jc$.

Suppose for contradiction that $W\in\Com\Jc$.
Let $\{\lambda_k\}_{k=1}^\infty$ be a rearrangement of $\{w_k\}_{k=1}^\infty$
such that $|\lambda_1|\ge|\lambda_2|\ge\cdots$.
Then by~\cite{1}, $\diag\{\frac1k(\lambda_1+\cdots+\lambda_k)\}\in\Jc$.
For every $k\in\Nats$ we have
$$ \left|\sum_{j=1}^k\lambda_j-\sum_{j=1}^kw_j\right|\le2ku_k $$
because every $w_j$ having absolute value strictly greater than $u_k$ appears
in both summations above.
Hence $\diag\{\frac1k(w_1+\cdots w_k)\}\in\Jc$, thus by~\scite{5}{3.1(3)} there
are $\alpha>0$ and $c>0$ such that
$$ \forall k\in\Nats\qquad\frac1k\theta_k\le cu_{\alpha k}. \tag{\thetasinth} $$
Let $k=2^{p_n+n}$.
We will find a lower bound for $\theta_k$ which will contradict~(\thetasinth).
Routine estimation reveals that
$$ jv_j+|\sigma_k-\sigma_j|\ge\cases
\frac12&\text{if }j<2^{p_{n-1}} \\
1&\text{if }2^{p_{n-1}}\le j\le2^{q_{n-1}}-1 \\
2^{-p_{n+1}+p_n+2n}&\text{if }2^{q_{n-1}}\le j\le2^{p_n}-1 \\
2^{-p_{n+1}+p_n+2n}&\text{if }2^{p_n}\le j\le2^{p_n+n}-1 \\
2^{-p_{n+1}+p_n+3n}&\text{if }2^{p_n+n}\le j\le2^{p_n+2n}-1 \\
2^{-p_{n+1}+p_n+2n-1}&\text{if }2^{q_n}\le j\le2^{p_{n+1}}-1 \\
\frac12&\text{if }2^{p_{n+1}}\le j.
\endcases $$
Considering all cases, we find
$\theta_k\ge2^{-p_{n+1}+p_n+2n-1}$, so $\frac1k\theta_k\ge2^{-p_{n+1}+n-1}$.
Let $\alpha>0$ be arbitrary.
If $n$ is so large that $2^n>\alpha^{-1}$ then
$u_{\alpha k}=2^{-p_{n+1}}$ and
$\frac1k\theta_k/u_{\alpha k}\ge2^{n-1}$, which grows without bound as
$n\to\infty$,
contradicting~(\thetasinth).
\QED

\heading \SpectralTraces.  Spectral traces. \endheading

Let $\Jc$ be an ideal of compact operators.
A {\it trace} on $\Jc$ is a linear functional, $\tau:\Jc\to\Cpx$ that is
unitarily invariant, or, equivalently,  that vanishes
on $\Com\Jc$.
In this section, we show that, as a consequence
of~\scite{5}{Theorem 3.3}, every trace on a geometrically stable ideal is
a spectral trace, i.e\. for every $T\in\Jc$, $\tau(T)$ depends only on the
eigenvalues of $T$, listed according to algebraic multiplicity.

In the following, given a (compact) operator, $T$, on some Hilbert space
$\Hil$, the spectrum of $T$ will be denoted $\sigma(T)$.
If $\lambda\in\sigma(T)\backslash\{0\}$ then $E_\lambda(T)$ will denote the
finite dimensional subspace, $\bigcup_{n\ge1}\ker(\lambda-T)^n$.
Note that $E_\lambda(T)$ is invariant under $T$, and has dimension equal to the
algebraic multiplicity of $\lambda$ as an eigenvalue of $T$.
We will denote by $P_\lambda$ the  orthogonal projection of $\Hil$ onto
$E_\lambda(T)$.

\proclaim{Lemma \SpecOut}
Let $T$ be a compact operator on infinite dimensional Hilbert space $\Hil$, and
suppose $\lambda\in\sigma(T)\backslash\{0\}$.
Then
$$ \sigma\bigl((1-P_\lambda)T(1-P_\lambda)\bigr)
=\sigma(T)\backslash\{\lambda\}. \tag{\specout} $$
\endproclaim
\demo{Proof}
Since both $T$ and its compression are compact operators on infinite
dimensional Hilbert space, $0$ is an element of both sides of~(\specout).
Suppose $\mu\in\sigma(T)\backslash\{0,\lambda\}$.
Then there is nonzero $\zeta\in\Hil$ such that $T\zeta=\mu\zeta$.
One easily sees that $(1-P_\lambda)\zeta$ is nonzero and is an eigenvector of
$(1-P_\lambda)T(1-P_\lambda)$ with eigenvalue $\mu$.
This shows that the inclusion $\supseteq$ holds in~(\specout).

To show the opposite inclusion, suppose
$\mu\in\sigma((1-P_\lambda)T(1-P_\lambda))\backslash\{0\}$.
Then there is nonzero $\zeta\in(1-P_\lambda)\Hil$ such that
$(1-P_\lambda)T\zeta=\mu\zeta$.
Thus
$$ \eta\eqdef(\mu-T)\zeta\in E_\lambda(T)T. \tag{\etalambda} $$
Suppose $\mu\ne\lambda$.
Since $P_\lambda TP_\lambda$ has spectrum $\{\lambda\}$, there is
$\eta'\in E_\lambda(T)T$ such that $(\mu-T)\eta'=\eta$.
Let $\xi=\zeta-\eta'$.
Then $\xi\ne0$ and $(\mu-T)\xi=0$;
hence $\mu\in\sigma(T)$.

If, on the other hand, $\mu=\lambda$, then from~(\etalambda) we have
$\zeta\in E_\lambda(T)T$, a contradiction.
\QED

Let $P$ be the projection of $\Hil$ onto
$$ \lspan\bigcup_{\lambda\in\sigma(T)\backslash\{0\}}E_\lambda(T). $$
Note that $PTP=TP$.
Choosing bases for all the $E_\lambda(T)$ and applying Gramm--Schmidt, we can find
an orthonormal basis for $P\Hil$ with respect to which $PTP$ is upper
triangular and has on the main diagonal all the elements of
$\sigma(T)\backslash\{0\}$ repeated according to their algebraic multiplicities.
Thus $\sigma(PTP)=\sigma(T)$.

\proclaim{Lemma \WhatsLeftNil}
Let $T$ be a compact operator on infinite dimensional Hilbert space $\Hil$,
and let $P$ be the projection defined above.
Then $(1-P)T(1-P)$ is quasi--nilpotent.
\endproclaim
\demo{Proof}
Suppose for contradiction there is nonzero $\mu\in\sigma((1-P)T(1-P))$.
Then there is nonzero $\zeta\in(1-P)\Hil$ such that
$\eta\eqdef(\mu-T)\zeta\in P\Hil$.
If $\mu\notin\sigma(T)$ then since $\sigma(T)=\sigma(PTP)$, there is $\eta'\in
P\Hil$ such that $(\mu-T)\eta'=\eta$.
Thus $\zeta-\eta'\ne0$ and $(\mu-T)(\zeta-\eta')=0$, contradicting that
$\mu\notin\sigma(T)$.

If, on the other hand, $\mu\in\sigma(T)$ then $\eta=\eta_\mu+\eta_o$ where
$\eta_\mu\in E_\mu(T)$ and $\eta_o\perp E_\mu(T)$.
By Lemma~\SpecOut, $\mu$ is not in the spectrum of $(P-P_\mu)T(P-P_\mu)$, so
there is $\eta'\in(P-P_\mu)\Hil$ such that $(P-P_\mu)(\mu-T)\eta'=\eta_o$.
Thus $(\mu-T)(\zeta-\eta')\in E_\mu(T)$, so $\zeta-\eta'\in E_\mu(T)$, which
contradicts the choice of $\zeta$.
\QED

\proclaim{Lemma \QuasiNTriang}
Let $T\in B(\Hil)$.
Suppose that for a projection $P$, $PTP$ and $(1-P)T(1-P)$ are
quasi--nilpotent and $(1-P)TP=0$.
Then $T$ is quasi--nilpotent.
\endproclaim
\demo{Proof}
For convenience write $P_1=P$ and $P_2=1-P$.
Then
$$ T^n=(P_1TP_1)^n+\sum_{k=0}^{n-1}(P_1TP_1)^{n-k-1}(P_1TP_2)(P_2TP_2)^k
+(P_2TP_2)^n. $$
Now some elementary estimates and the quasi--nilpotence of $P_1TP_1$ and
$P_2TP_2$  show that $T$ is quasi--nilpotent.
\QED

Putting together these lemmas, we have the following.
\proclaim{Proposition \DiagQnil}
Let $T$ be a compact operator on infinite dimensional Hilbert space $\Hil$.
Then $T=D+Q$, where $D$ is a normal operator whose eigenvalues and
multiplicities are equal to those of $T$, and where $Q$ is a quasi--nilpotent
operator.
\endproclaim

\proclaim{Corollary \GSDiagQnil}
Let $\Jc$ be a geometrically stable ideal of compact operators and let
$T\in\Jc$.
Then $T=D+Q$ where $D\in\Jc$ is normal and $Q\in\Jc$ is quasi--nilpotent.
\endproclaim
\demo{Proof}
If $\lambda_k=\lambda_k(T)$ are the eigenvalues of $T$ listed according to
algebraic multiplicity, then the inequality
$|\lambda_1\cdots\lambda_k|\le|s_1(T)\cdots s_k(T)|$ and the geometric
stability of $\Jc$ imply that $D\in\Jc$.
\QED

\proclaim{Corollary \SpecTrace}
Let $\Jc$ be a geometrically stable ideal and suppose $\tau$ is a trace on
$\Jc$.
For every given $T\in\Jc$, the value $\tau(T)$ depends only on the eigenvalues
of $T$ and their algebraic multiplicities.
\endproclaim
\demo{Proof}
Using the decomposition $T=D+Q$ from Corollary~\GSDiagQnil{} and the spectral
characterization of $\Com\Jc$ in~\scite{5}{3.3}, it follows that
$Q\in\Com\Jc\subseteq\ker\tau$.
Hence $\tau(T)=\tau(D)$.
\QED

\enddocument